\input amstex
\magnification=\magstep1
\documentstyle{amsppt}
\pagewidth{5.5 true in}
\pageheight{8.0 true in}
\leftheadtext{Y. J. Cho, S. S. Dragomir and  Y. H. Kim}
\rightheadtext{On some  Gronwall type  inequalities}
\topmatter
\title 
 On some  Gronwall type  inequalities involving iterated integrals
\endtitle
\endtopmatter
\centerline{Yeol Je Cho}
\centerline{Department of Mathematics Education}
\centerline{The Research Institute of Natural Sciences}
\centerline{College of Education, Gyeongsang National University}
\centerline{Chinju 660-701, Republic of Korea}
\centerline{\it E-mail: yjcho\@nongde.gsnu.ac.kr}
\vskip 2mm

\centerline{Sever S. Dragomir}
\centerline{School of Computer Science and Mathematics}
\centerline{Victoria University of Technology}
\centerline{PO Box 14428, MCMC, Melbourne}
\centerline{Victoria 8001, Australia}
\centerline{\it E-mail: sever.dragomir\@vu.edu.au}
\vskip 2mm

\centerline{Young-Ho Kim}
\centerline{Department of Applied  Mathematics}
\centerline{Changwon National University}
\centerline{Changwon 641-773, Republic of Korea}
\centerline{\it E-mail:  yhkim\@sarim.changwon.ac.kr}
\vskip 5mm

\noindent
{\bf Abstract:} In this paper,  some new Gronwall type inequalities involving iterated integrals are given.
\vskip 2mm

\footnote""{{\bf 2000 AMS Subject Classification:}  26D15, 35A05.}
\footnote""{{\bf Key Words and Phrases:} Gronwall-Bellman inequalities, integral inequality, iterated integrals, nondecreasing function.}
\document
\NoBlackBoxes
\bigskip

\centerline{\bf 1. Introduction}
\vskip 2mm

Let $u:[\alpha , \alpha +h] \to R$ be a continuous real-valued function satisfying the inequality
$$
0\leq u(t) \leq \int_{\alpha}^t [a + bu(s)] \, ds,  \quad   t\in [\alpha , \alpha +h],
$$
where $a,b$ are nonnegative constants. Then $u(t)\leq ahe^{bh}$ for $t\in [\alpha, \alpha +h].$ This result was proved by T. H.
Gronwall [9] in the year 1919, and is the prototype for the study of several integral inequalities of Volterra type, and also for obtaining
explicit bounds of the unknown function. Among the several results on this subject, the inequality of Bellman [3] is very well known:
\vskip 2mm

{\it Let $x(t)$ and $k(t)$ be real valued nonnegative continuous functions for $t\geq \alpha .$ If $a$ is a constant, $a\geq 0,$ and
$$
x(t)\leq a + \int_{\alpha}^t
 k(s)x(s)  \, ds,\quad   t \geq \alpha,
$$
then
$$
x(t)\leq a \exp\biggl( \int_{\alpha}^t
 k(s)   \, ds\biggr) ,\quad   t \geq \alpha .
$$}
\vskip 2mm

It is clear that Bellman's result contains that of Gronwall. This is the reason why inequalities of this type were called ``Gronwall-Bellman
inequalities" or ``Inequalities of Gronwall type". The Gronwall type integral inequalities provide a necessary tool for the study of the theory of differential equations, integral equations and inequalities of various types (see Gronwall [9] and Guiliano [10]). Some applications of this result to the study of stability of
the solution of linear and nonlinear differential equations may be found in Bellman [3]. Some applications to existence and uniqueness theory of
differential equations may be found in Nemyckii-Stepanov [14], Bihari [4], and Langenhop [11]. During the past few years several authors (see references below
and some of the references cited therein) have established several Gronwall type integral inequalities in two or more independent real variables. Of course, such results have application in the theory of partial differential equations and Volterra integral equations.

Bykov proved the following interesting integral inequality, which appear in [1, p. 98]:
\vskip 2mm

{\it Let $u(t),$  $ b(t),$ $ k(t,s)$ and $h(t,s, \sigma )$ be nonnegative continuous functions for $\alpha\leq \tau \leq s\leq t\leq \beta$ and suppose that
$$
\aligned
u(t) & \leq a +\int_\alpha^t b(s )u(s )\, ds + \int_\alpha^t \int_\alpha^s k(s,\tau ) u(\tau) \, d\tau  ds  \\
&\quad + \int_\alpha^t \int_\alpha^s \int_\alpha^\tau h(s,\tau , \sigma ) u(\sigma) \, d\sigma \, d\tau  ds
\endaligned
 \tag{1.1}
$$
for any $t \in [\alpha, \beta],$ where $a\geq 0$ is a constant. Then
$$
\aligned
u(t) &\leq a \exp\biggl( \int_\alpha^t b(s ) \, ds + \int_\alpha^t \int_\alpha^s k(s,\tau ) \, d\tau  ds \\
&\quad  + \int_\alpha^t \int_\alpha^s \int_\alpha^\tau h(s,\tau , \sigma)\, d\sigma \, d\tau \,ds \biggr),\quad t \in [\alpha, \beta].
\endaligned
$$
}
\vskip 2mm

In this paper,  we consider simple inequalities involving iterated integrals in the inequality (1.1) for the case when the function $u$
in the right-hand side of the inequality (1.1) is replaced by the function $u^p$ for some $p,$ and  the constant $a$ is replaced by a nonnegative, nondecreasing
function $a(t)$. We also provide some related integral inequalities involving iterated integrals.
\bigskip

\centerline{\bf 2. The case  $p>1$}
\vskip 2mm

In this section, we state and prove some new nonlinear integral inequalities  involving iterated integrals. Throughout the paper,
all the functions which appear in the inequalities are assumed to be real-valued. Before considering our first integral inequality involving iterated
integrals, we need the following lemma, which appears in [1, p. 2].
\vskip 2mm

\proclaim{ Lemma 2.1} Let $b(t)$ and $f(t)$ be continuous function for $t\geq \alpha,$ let $v(t)$ be a differentiable function for $t\geq \alpha$ and
suppose that
$$
v'(t)\leq b(t)v(t) +f(t), \quad t\geq \alpha,
$$
and $v(\alpha ) \leq v_0.$ Then we have
$$
v(t)\leq v_0 \exp \biggl(\int_{\alpha }^{t}b(s) \, ds\biggr)
+\int_{\alpha }^{t}f(s)\exp \biggl(\int_{s }^{t} b(\tau )\,
d\tau\biggr)\, ds, \quad t\geq \alpha.
$$
\endproclaim
\vskip 2mm

\proclaim{ Theorem 2.2} Let $u(t),$  $ b(t),$ $ k(t,s)$ and $h(t,s, \tau )$ be nonnegative continuous functions for $\alpha\leq \tau\leq s\leq t\leq \beta$ and let
$p>1$ be a constant. Suppose $a(t)\geq 0$ is nondecreasing in $J=[\alpha, \beta ]$ and
$$
\aligned
u(t) &\leq a(t)+\int_\alpha^t b(s )u^p(s )\, ds + \int_\alpha^t \int_\alpha^s k(s,\tau ) u^p(\tau) \, d\tau  ds   \\
&\quad + \int_\alpha^t \int_\alpha^s \int_\alpha^\tau h(s,\tau , \sigma ) u^p(\sigma) \, d\sigma \, d\tau  ds,\quad t \in [\alpha, \beta].  
\endaligned
\tag{2.1}
$$
Then we have
$$
u(t) \leq a(t) \biggl[1-(p-1)\int_\alpha^t B(s) a^{p-1}(s) \, ds
 \biggr]^{1\over 1-p},\quad  t \in [\alpha, \beta_p ), 
\tag{2.2}
$$
where
$$
\beta_p = \sup\Big\{t\in J:(p-1)\int_\alpha^t B(s) a^{p-1}(s) \,ds<1\Big\}
$$
and
$$
B(t)= b(t ) + \int_\alpha^t  k(t,s ) \, ds   + \int_\alpha^t \int_\alpha^s  h(t,s , \tau ) \, d\tau \,ds. 
$$
\endproclaim
\vskip 2mm

\demo{ Proof } We denote the right-hand side of (2.1) by  $a(t) + v(t).$ Then, for $\alpha\leq t\leq T<\beta_p,$ (2.1) implies $v(\alpha )=0,$
the function $v(t)$ is nondecreasing in  $t \in [\alpha, \beta],$
$$
u(t)\leq a(t) + v(t)
\tag{2.3}
$$
and
$$
\aligned
 v' (t)
&= b(t )u^p(t )  + \int_\alpha^t   k(t,\tau )
 u^p(\tau) \, d\tau   + \int_\alpha^t   \int_\alpha^\tau h(t,\tau , \sigma )
 u^p(\sigma) \, d\sigma \, d\tau     \\
&\leq B(t )[a(t) + v(t )  ]^p     \\
&\leq B(t )[a(t) + v(t )  ]^{p-1}[a(T) + v(t )  ],\\
\endaligned
$$
that is,
$$
 v' (t)\leq R(t )[a(T) + v(t )  ],  
\tag{2.4}
$$
where $R(t)=B(t )[a(t) + v(t )  ]^{p-1}.$ Lemma 2.1 and (2.4) imply
$$
v(t)\leq a(T)\int_{\alpha}^t R(s) \exp \biggl(\int_s^tR(\tau )\,
d\tau\biggr)\, ds
$$
and so
$$
v(t)+a(T)\leq a(T) \exp \biggl(\int_\alpha^tR(s )\,
ds\biggr), \quad \alpha \leq t\leq T.
$$
Hence, for $t=T,$
$$
v(t)+a(t)\leq a(t) \exp \biggl(\int_\alpha^tR(s )\,
ds\biggr) .
\tag{2.5}
$$
From (2.5),  we successively obtain
$$
\aligned
&[v(t)+a(t)]^{p-1}\leq a^{p-1}(t) \exp \biggl(\int_\alpha^t(p-1)R(s )\,ds\biggr),\\
&R(t)\leq B(t)a^{p-1}(t) \exp \biggl(\int_\alpha^t(p-1)R(s )\,ds\biggr),\\
&Z(t) \leq (p-1)B(t)a^{p-1}(t) \exp \biggl(\int_\alpha^t Z(s )\,ds\biggr),
\endaligned
$$
where $Z(t)=(p-1)R(t).$ Consequently, we have
$$
Z(t)\exp \biggl(- \int_\alpha^t Z(s )\,ds\biggr)\leq (p-1)B(t)a^{p-1}(t)
$$
or
$$
\frac{d}{dt} \biggl[- \exp \biggl(- \int_\alpha^t Z(s )\,ds\biggr)\biggr]\leq (p-1)B(t)a^{p-1}(t).
$$
Integrating this from $\alpha $ to $t$ yields
$$
1- \exp \biggl(- \int_\alpha^t Z(s )\, ds\biggr)\leq \int_{\alpha}^{t} (p-1)B(s)a^{p-1}(s) \, ds,
$$
from which we conclude that
$$
  \exp \biggl(  \int_\alpha^t R(s )\, ds\biggr)
\leq \biggl[1-(p-1)\int_{\alpha}^{t}  B(s)a^{p-1}(s) \,
ds\biggr ]^{\frac{1}{1-p}}.
$$
This, together with (2.3) and (2.5), implies (2.2). This completes the proof.
\enddemo
\vskip 2mm

In the same manner,  we can prove the following theorem:
\vskip 2mm

\proclaim{ Theorem 2.3} Let $u(t),$  $ b(t),$ $ k(t,s)$ and $\sigma(t)$ be nonnegative continuous functions for $\alpha\leq s\leq t\leq \beta$ and let
$p>1$ be a constant. Suppose that $\sigma (t)$ is nondecreasing in $J=[\alpha, \beta ]$ and
$$
u(t) \leq \sigma (t)\biggl\{a_1 +\int_\alpha^t b(s )u^p(s )\, ds + \int_\alpha^t \int_\alpha^s k(s,\tau ) u^p(\tau) \, d\tau  ds\biggr\}
$$
for any $t \in [\alpha, \beta],$ where $a_1\geq 0$ is a constant. Then we have
$$
u(t) \leq a_1\sigma (t)\exp (\sigma (t))  \biggl[1-(p-1)a_1^{p-1}
\int_\alpha^t B_1(s) \sigma^{p-1}  (s)\exp ( \sigma (s)) \, ds \biggr]^{1\over 1-p}
$$
for any $t \in [\alpha, \beta_p ),$ where
$ B_1(t)= b(t ) + \int_\alpha^t  k(t,\tau ) \, d\tau $ and
$$
 \beta_p = \sup\{t\in J:(p-1)a_1^{p-1}\int_\alpha^t B_1 (s) \sigma^{p-1}(s)\exp ( \sigma (s)) \,ds<1\}.
$$
\endproclaim
\vskip 2mm

Let $\alpha <\beta ,$ and set $J_i =\{(t_1, t_2, \ldots ,t_i )\in R^i:\alpha\leq t_i\leq \cdots \leq t_1\leq \beta \},$ $ i=1,\cdots , n.$
\vskip 2mm

\proclaim{Theorem 2.4} Let $u(t), a(t)$ and  $   b(t) $ be nonnegative continuous functions in $J =[\alpha , \beta ]$ and let $p>1$ be a constant. Suppose that $\frac{a(t)}{b(t)}$ is nondecreasing in $J$ and
$$
\aligned
u(t)& \leq a(t)+ b(t) \biggl[ \int_{\alpha}^t k_1(t,t_1)u^p(t_1)\, dt_1 + \cdots \\
&\quad + \int_{\alpha}^t\biggl(\int_{\alpha}^{t_1}\cdots\biggl( \int_{\alpha}^{t_{n-1}}k_n(t,t_1, \cdots ,t_n)u^p(t_n)  \, dt_n\biggr) \cdots\biggr) dt_1\biggr]\\
\endaligned
\tag{2.6}
$$
for any $t\in J,$ where $k_i(t, t_1,\ldots ,t_i) $ are nonnegative continuous functions in $J_{i+1}$ for $i=1,2,\cdots,n.$ Suppose thta the partial derivatives $\frac{\partial k_i}{\partial t}(t, t_1,\cdots,t_i)$ exist and are nonnegative and continuous in $J_{i+1}$ for $i=1,2,\cdots ,n.$ Then, for any $t\in J,$
$$
u(t) \leq a(t)\biggl[1-(p-1) \int_{\alpha }^t \biggl(\frac{a(s)}{b(s)}\biggr)^{p-1}  (R[b^p] (s)+Q[b^p](s)) \, ds\biggr]^{ \frac{1}{1-p}} 
\tag{2.7}
$$
for any $t \in [\alpha, \tilde{\beta_p}  ),$ where
$$
  \tilde{\beta_p}= \sup \{t\in J:(p-1)a_1^{p-1}\int_\alpha^t( {a(s)}/{b(s)} )^{p-1}  (R[b^p] (s)+Q[b^p](s) \,ds<1  \},
$$
$$
\aligned
  R[w](t)& =k_1(t,t)w(t) +\int_{\alpha }^t k_2(t,t,t_2) w(t_2) dt_2 \\
&\quad + \sum_{i=3}^{n}\int_{\alpha }^t\biggl(\int_{\alpha }^{t_2} \cdots \biggl(\int_{\alpha }^{t_{i-1}}k_i(t,t, t_2, \cdots ,t_i)w(t_i) \, dt_i\biggr)
 \cdots \biggr) \, dt_2,\\
 Q[w](t)& =\int_{\alpha }^t \frac{\partial k_1}{\partial t} (t,t_1)w(t_1)\, dt_1 \\
&\quad + \sum_{i=2}^{n}\int_{\alpha }^t \biggl(\int_{\alpha }^{t_1} \cdots \biggl(\int_{\alpha }^{t_{i-1}}\frac{\partial  k_i}{\partial t}(t,t_1,  \cdots ,t_i)w(t_i) \,
 dt_i\biggr) \cdots \biggr) \, dt_1
\endaligned
$$
for each continuous function $w(t)$ in $J.$
\endproclaim
\vskip 2mm

\demo{ Proof } First, we note that $R[w]$ and $Q[w]$ are linear functionals and
$$
R[w_1]\leq R[w_2], \quad Q[w_1]\leq Q[w_2]
$$
if $w_1(t) \leq w_2(t)$ for any $t\in J$ and
$$
R[w_1 w_2]\leq R[w_1]w_2, \quad Q[w_1w_2]\leq Q[w_1]w_2
$$
if $w_1(t)$ is nonnegative in $J$ and $w_2(t)$ is nondecreasing in
$J.$
We set
$$
\aligned
v(t)&= \int_{\alpha}^t k_1(t,t_1)u^p(t_1)\, dt_1 + \cdots \\
&\quad + \int_{\alpha}^t\biggl(\int_{\alpha}^{t_1}\cdots\biggl( \int_{\alpha}^{t_{n-1}}k_n(t,t_1, \cdots ,t_n)u^p(t_n)  \, dt_n\biggr)
 \cdots\biggr) dt_1.
\endaligned
$$
Then, for $\alpha\leq t\leq T<\beta_p,$ (2.6) implies $v(\alpha )=0,$ the function $v(t)$ is nondecreasing,
$$
u(t)\leq a(t) + b(t)v(t)
\tag{2.8}
$$
and we have
$$
v'(t)=R[u^p](t) +Q[u^p](t) \leq (R[b^p ](t) +Q[b^p ](t))(\frac{a(t)}{b(t)} +v(t))^p,
$$
that is,
$$
 v' (t)\leq R(t )[a(T)/b(T) + v(t )  ],  
\tag{2.9}
$$
where $R(t)=(R[b^p ](t) +Q[b^p ](t))[a(t)/b(t) + v(t )  ]^{p-1}.$ Lemma 2.1 and (2.9) imply
$$
v(t)+\frac{a(T)}{b(T)}\leq \frac{a(T)}{b(T)} \exp \biggl(\int_\alpha^tR(s )\,ds\biggr), \quad \alpha \leq t\leq T.
$$
Hence, for $t=T,$
$$
v(t)+\frac{a(t)}{b(t)}\leq \frac{a(t)}{b(t)} \exp \biggl(\int_\alpha^tR(s )\,ds\biggr).
\tag{2.10}
$$
From (2.10), we successively obtain
$$
\aligned
&\biggl[v(t)+\frac{a(t)}{b(t)}\biggr]^{p-1}\leq\biggl[\frac{a(t)}{b(t)}\biggr]^{p-1}  \exp \biggl(\int_\alpha^t(p-1)R(s )\,ds\biggr),\\
&R(t)\leq (R[b^p ](t) +Q[b^p ](t)) \biggl[\frac{a(t)}{b(t)}\biggr]^{p-1}\exp \biggl(\int_\alpha^t(p-1)R(s )\,ds\biggr),\\
&Z(t) \leq (p-1)(R[b^p ](t) +Q[b^p ](t)) \biggl[\frac{a(t)}{b(t)}\biggr]^{p-1}\exp \biggl(\int_\alpha^t(p-1)R(s )\,ds\biggr),
\endaligned
$$
where $Z(t)=(p-1)R(t).$ Consequently, we have
$$
\frac{d}{dt} \biggl[- \exp \biggl(- \int_\alpha^t Z(s )\,ds\biggr)\biggr]
\leq (p-1)(R[b^p ](t) +Q[b^p ](t)) \biggl[\frac{a(t)}{b(t)}\biggr]^{p-1}.
$$
Integrating this from $\alpha $ to $t$ yields
$$
\aligned
&1- \exp \biggl(- \int_\alpha^t Z(s )\, ds\biggr)\\
&\leq (p-1) \int_{\alpha}^{t} \biggl(\frac{a(s)}{b(s)}\biggr)^{p-1}(R[b^p ](s) +Q[b^p ](s)) \, ds,
\endaligned
$$
from which we conclude that
$$
\aligned
&  \exp \biggl(  \int_\alpha^t R(s )\, ds\biggr)\\
&\leq \biggl[1-(p-1)\int_{\alpha}^{t} \biggl(\frac{a(s)}{b(s)}\biggr)^{p-1}(R[b^p ](s) +Q[b^p ](s)) \,ds\biggr ]^{\frac{1}{1-p}}.
\endaligned
$$
This, together with (2.8) and (2.10), implies (2.7). This completes the proof.
\enddemo
\bigskip

\centerline{\bf 3. The case  $p>0$ $(p\neq 1)$}
\vskip 2mm

In this section, we use another method for studying nonlinear integral inequalities. Before considering the first result of the integral inequality, we
need the following lemma, which appears in [1, p. 38].
\vskip 2mm

\proclaim{ Lemma 3.1} Let $v(t)$ be a positive differential function satisfying the inequality
$$
v'(t)\leq b(t)v(t) + k(t)v^p(t), \quad t\in J=[\alpha, \beta ],
$$
where the functions $b$ and $k$ are continuous in $J$ and $p\geq 0$ ($p\neq 1$) is a constant. Then we have
$$
v(t)\leq \exp\biggl(\int_{\alpha}^{t}b(s)\, ds\biggr)\biggl[v^q(\alpha ) + q\int_{\alpha}^{t}k(s)\exp\biggl(-q\int_{\alpha}^{s}b(\tau)\,
d\tau\biggr)\, ds\biggr]^{1/q}
$$
for any $t\in [\alpha ,\beta_1 ),$ where $\beta_1$ is chosen so that the expression between $[\cdots]$ is
positive in the subinterval $[\alpha , \beta_1 ).$
\endproclaim
\vskip 2mm

An essential element in the investigation of the integral inequalities in the following theorems is the application of the result of Lemma 3.1.
\vskip 2mm

\proclaim{ Theorem 3.2} Let $u(t),$  $ b(t),$ $ k(t,s),$ $h(t,s, \sigma )$ be nonnegative continuous functions for $\alpha\leq \sigma\leq s\leq t\leq \beta$ and suppose that
$$
\aligned
u(t) &\leq a+\int_\alpha^t b(s )u^p(s )\, ds + \int_\alpha^t \int_\alpha^s k(s,\tau ) u^p(\tau) \, d\tau  ds   \\
&\quad + \int_\alpha^t \int_\alpha^s \int_\alpha^\tau h(s,\tau , \sigma ) u^p(\sigma) \, d\sigma \, d\tau  ds
\endaligned
\tag{3.1}
$$
for any $t \in [\alpha, \beta],$ where $a> 0$ and $p\geq 0$ $(p\neq 1)$ are a constants.
Then we have
$$
\aligned
u(t) &\leq\biggl[a^q +q \int_\alpha^t \biggl( b(s )  +  \int_\alpha^s k(s,\tau )   \, d\tau   \\
&\quad   +   \int_\alpha^s \int_\alpha^\tau h(s,\tau , \sigma )  \, d\sigma \, d\tau \biggr)\, ds \biggr]^{1/q} 
\endaligned
 \tag{3.2}
$$
for any $t \in [\alpha, \beta_1 ),$ where $q=1-p$ and $\beta_1$ is chosen so that the expression between $[\cdots]$ is positive in the subinterval 
$ [\alpha , \beta_1).$
\endproclaim
\vskip 2mm

\demo{ Proof } We denote the right-hand side of (3.1) by the function $v(t).$ Then the function $v(t)$ is nondecreasing in  $t \in [\alpha, \beta],$
$u(t)\leq v(t),$ $v(\alpha )=a$ and
$$
\aligned
 v' (t)&= b(t )u^p(t )  + \int_\alpha^t   k(t,\tau ) u^p(\tau) \, d\tau   + \int_\alpha^t   \int_\alpha^\tau h(t,\tau , \sigma ) u^p(\sigma) \, d\sigma \, d\tau     \\
&\leq b(t )v^p(t )  + \int_\alpha^t   k(t,\tau ) v^p(\tau) \, d\tau   + \int_\alpha^t   \int_\alpha^\tau h(t,\tau , \sigma ) v^p(\sigma) \, d\sigma \, d\tau      \\
&\leq \biggl(b(t )  + \int_\alpha^t   k(t,\tau ) \, d\tau   + \int_\alpha^t   \int_\alpha^\tau h(t,\tau , \sigma ) \, d\sigma \, d\tau \biggr) v^p(t ).  
\endaligned
$$
Therefore, applying Lemma 3.1,  we arrive at (3.2). This completes the proof.
\enddemo
\vskip 2mm

\proclaim{Theorem 3.3} Let $u(t),$  $ b(t),$ $ k(t,s)$, $h(t,s, \sigma )$ be nonnegative continuous functions for $\alpha\leq \sigma\leq s\leq t\leq \beta$ and suppose that
$$
\aligned
u(t)& \leq a(t)+\int_\alpha^t b(s )u^p(s )\, ds + \int_\alpha^t \int_\alpha^s k(s,\tau ) u^p(\tau) \, d\tau  ds   \\
&\quad + \int_\alpha^t \int_\alpha^s \int_\alpha^\tau h(s,\tau , \sigma ) u^p(\sigma) \, d\sigma \, d\tau  ds
\endaligned
\tag {3.3}
$$
for any $t \in [\alpha, \beta],$ where $a(t)$ is a positive nondecreasing function and $p\geq 0$ $(p\neq 1)$ is a constant.
Then we have
$$
\aligned
u(t) &\leq\biggl[A^q(t) +q \int_\alpha^t \biggl( b(s )  +  \int_\alpha^s k(s,\tau )   \, d\tau  \\
&\quad    +   \int_\alpha^s \int_\alpha^\tau h(s,\tau , \sigma )  \, d\sigma \, d\tau \biggr)\, ds \biggr]^{1\over q}
\endaligned
  \tag {3.4}
$$
for any $t \in [\alpha, \beta_1 ),$ where $q=1-p,$ $A(t)=\sup_{s\in [\alpha, t]} a(s)$ and $\beta_1$ is chosen so that the expression between
$[\cdots]$ is positive in the subinterval $  [\alpha , \beta_1) .$
\endproclaim
\vskip 2mm

\demo{ Proof } The function $A(t)$ is nondecreasing in $t \in [\alpha, \beta ].$ Thus (3.3) implies that, for all $\alpha\leq t\leq T\leq \beta ,$
$$
\aligned
u(t) &\leq A(T)+\int_\alpha^t b(s )u^p(s )\, ds + \int_\alpha^t \int_\alpha^s k(s,\tau ) u^p(\tau) \, d\tau  ds    \\
&\quad + \int_\alpha^t \int_\alpha^s \int_\alpha^\tau h(s,\tau , \sigma ) u^p(\sigma) \, d\sigma \, d\tau  ds. 
\endaligned
\tag{3.5}
$$
We denote the right-hand side of (3.5) by the function $v(t).$ Then the function $v(t)$ is nondecreasing in  $t \in [\alpha, \beta],$ $u(t)\leq v(t),$
$v(\alpha )=A(T)$
and
$$
v' (t)\leq \biggl(b(t )  + \int_\alpha^t   k(t,\tau ) \, d\tau   + \int_\alpha^t   \int_\alpha^\tau h(t,\tau , \sigma ) \, d\sigma \, d\tau \biggr) v^p(t ).
$$
Consequently, Lemma 3.1 implies
$$
u(t) \leq\biggl[A^q(T) +q\int_\alpha^t \biggl( b(s )  +  \int_\alpha^s k(s,\tau )   \, d\tau      +   \int_\alpha^s \int_\alpha^\tau h(s,\tau , \sigma )
  \, d\sigma \, d\tau \biggr)\, ds
 \biggr]^{1\over q}
$$
and, for  $t=T$, we obtain (3.4). This completes the proof.
\enddemo
\vskip 2mm

Let $\alpha <\beta ,$ and set 
$$
J_i =\{(t_1, t_2, \ldots ,t_i )\in R^i :\alpha\leq t_i\leq \cdots \leq t_1\leq \beta \}
$$ 
for $ i=1,\cdots , n.$ By a similar  reasoning to the proof of Theorem 3.2, we also can prove the following result:
\vskip 2mm

\proclaim{Theorem 3.4} Let $u(t),$ and  $   b(t) $ be nonnegative continuous functions in $J =[\alpha , \beta ]$ and suppose that
$$
\aligned
u(t) &\leq  b(t) \biggl[a+ \int_{\alpha}^t k_1(t,t_1)u^p(t_1)\, dt_1 + \cdots \\
&\quad + \int_{\alpha}^t\biggl(\int_{\alpha}^{t_1}\cdots\biggl( \int_{\alpha}^{t_{n-1}}k_n(t,t_1, \cdots ,t_n)u^p(t_n)  \, dt_n\biggr) \cdots
\biggr) dt_1\biggr]
\endaligned
$$
for any $ t\in J,$ where $a>  0$ and  $p\geq 0$ $(p\neq 1)$ is a constant, $k_i(t, t_1,\cdots ,t_i) $ are nonnegative continuous functions in
$J_{i+1}$ for $i=1,2,\cdots ,n.$ Suppose that the partial derivatives $\frac{\partial k_i}{\partial t}(t, t_1,\cdots,t_i)$ exist and are nonnegative and continuous in $J_{i+1}$ for $ i=1,2,\cdots ,n.$ Then, for any $t\in J$,
$$
u(t)\leq b(t)\biggl[ a^q+q \int_{\alpha }^t (R[b^p] (s)+Q[b^p](s)) \, ds\biggr]^{1/q}
 \tag {3.6}
$$
for any $t \in [\alpha, \beta_1 ),$ where $q=1-p,$ $\beta_1$ is chosen so that the expression between $[\cdots]$ is positive in the subinterval $[\alpha , \beta_1),$
$$
\aligned
R[w](t)&=k_1(t,t)w(t) +\int_{\alpha }^t k_2(t,t,t_2) w(t_2) dt_2 \\
&\quad + \sum_{i=3}^{n}\int_{\alpha }^t\biggl(\int_{\alpha }^{t_2} \cdots \biggl(\int_{\alpha }^{t_{i-1}}k_i(t,t, t_2, \cdots ,t_i)w(t_i) \, dt_i\biggr) \cdots
 \biggr) \, dt_2,\\
Q[w](t)&=\int_{\alpha }^t \frac{\partial k_1}{\partial t} (t,t_1)w(t_1)\, dt_1 \\
&\quad + \sum_{i=2}^{n}\int_{\alpha }^t  \biggl(\int_{\alpha }^{t_1} \cdots \biggl(\int_{\alpha }^{t_{i-1}}\frac{\partial  k_i}{\partial t}(t,t_1,  \cdots ,t_i)w(t_i) \,
 dt_i\biggr) \cdots \biggr) \, dt_1
\endaligned
$$
for each continuous function $w(t)$ in $J.$
\endproclaim
\vskip 2mm

\demo{ Proof } We set
$$
\aligned
v(t)&=a+ \int_{\alpha}^t k_1(t,t_1)u^p(t_1)\, dt_1 + \cdots \\
&\quad + \int_{\alpha}^t\biggl(\int_{\alpha}^{t_1}\cdots\biggl( \int_{\alpha}^{t_{n-1}}k_n(t,t_1, \cdots ,t_n)u^p(t_n)  \, dt_n\biggr) \cdots\biggr) dt_1.
\endaligned
$$
Since $v(\alpha )= a,$ $u(t)\leq b(t)v(t)$ and $v(t)$ is nondecreasing and continuous in $J,$ we have
$$
\aligned
v'(t)&=R[u^p](t) +Q[u^p](t) \leq R[b^pu^p](t) +Q[b^pu^p](t)\\
&\leq (R[b^p ](t) +Q[b^p ](t))v^p(t),
\endaligned
$$
from which, by the same method as in the proof of Theorem 3.2, we find the inequality (3.6). This completes the proof.
\enddemo
\vskip 2mm

\proclaim{Corollary 3.5} Let $u(t) $  be nonnegative  continuous function for $ \alpha \leq t \leq  \beta$and suppose that
$$
u(t) \leq a+ \int_{\alpha}^t k_1(t,s)u^p(s)\, ds  + \int_{\alpha}^t\biggl(\int_{\alpha}^{s}   h(t,s, \sigma )u^p(\sigma)  \, d\sigma\biggr)    ds,
$$
where $a> 0$ and  $p\geq 0$ $(p\neq 1)$ is a constant, $k(t, s) $ and $h(t, s, \sigma ) $ are nonnegative continuous functions for 
$\alpha \leq \sigma \leq s \leq t \leq  \beta .$ Suppose that the partial derivatives $\frac{\partial k }{\partial t}(t, s)$ and
 $\frac{\partial h }{\partial t}(t, s, \sigma)$  exist and are nonnegative and continuous for $\alpha \leq \sigma \leq s \leq t \leq  \beta .$
Then, for any $t\in J,$
$$
u(t)
\leq \biggl[ a^q+q \int_{\alpha }^t (R  (s)+Q (s)) \, ds\biggr]^{1/q},\quad t \in [\alpha, \beta_1 ),
$$
where $q=1-p,$ $\beta_1$ is chosen so that the expression between $[\cdots]$ is positive in the subinterval $  [\alpha , \beta_1) ,$
$$
R (t)=k (t,t)  +\int_{\alpha }^t h(t,t,\sigma)   d\sigma  
$$
and
$$
Q(t)=\int_{\alpha }^t \frac{\partial k}{\partial t} (t,\sigma )\, d\sigma+\int_{\alpha }^t
 \biggl(\int_{\alpha }^{s} \frac{\partial  h}{\partial t}(t,s, \sigma )  \, d\sigma \biggr)   \, ds .
$$
\endproclaim
\vskip 4mm

\centerline{\bf Acknowledgement}  
\vskip 2mm

Y. J. Cho and S. S. Dragomir greatly acknowledge the financial support from the Brain Pool Program (2002) of the Korean Federation of Science and  Technology Societies. The research was performed under the ``Memorandum of Understanding"  between Victoria University and Gyeongsang National University.

\enddocument